\documentclass[12pt,twoside]{article}

\usepackage{graphicx}
\usepackage{amssymb,amsmath}
\usepackage{epsfig}
\usepackage{epstopdf}

\newcommand{\1}{\mathds{1}}

\newcommand{\T}{{\mathcal T}}
\newcommand{\Tb}{{\mathbb T}}

\newcommand{\N}{{\mathbb N}}
\newcommand{\E}{{\mathbb E}}

\renewcommand{\P}{{\mathbb P}}

\renewcommand{\P}{{\mathbb P}}

\newcommand{\0}{{\mathbf 0}}

\newtheorem{stelling}{Theorem}

\numberwithin{equation}{section}
\allowdisplaybreaks

\begin{document}




\author{\footnotesize\textbf{Eric Cator, Henk Don}\\
\\
\footnotesize Institute for Mathematics, Astrophysics, and Particle Physics\\
\footnotesize Faculty of Science, Radboud University Nijmegen; e-mail:e.cator@math.ru.nl, henkdon@gmail.com}

\title{Conditioned multi-type Galton-Watson trees}
\maketitle
\begin{abstract}
We consider multi-type Galton Watson trees, and find the
distribution of these trees when conditioning on very general
types of recursive events. It turns out that the conditioned tree
is again a multi-type Galton Watson tree, possibly with more types
and with offspring distributions depending on the type of the
father node and on the height of the father node. These
distributions are given explicitly. We give some interesting
examples for the kind of conditioning we can handle, showing that
our methods have a wide range of applications.
\end{abstract}

\section{Introduction}

The asymptotic shape of conditioned Galton-Watson trees has been
widely studied. For example, one could condition on the number of
nodes of the tree being $n$, and letting $n\rightarrow\infty$. In
a lot of cases, the limiting tree is quite well understood, see
for example the survey paper by Janson \cite{Janson}. Some work on
finite conditioned trees has been done by Geiger and Kersting
\cite{Geiger}, who studies the shape of a tree conditioned on
having height exactly equal to $n$. In this context, we also
mention the spinal construction of a Galton-Watson tree
conditioned to reach generation k, as derived in \cite{Geiger2}.

In this paper, we will investigate conditioning multi-type
Galton-Watson trees on events of a recursive nature (as explained
in Section \ref{section:conditioning}), one example being
conditioning on survival to a given level. The main idea is that
we consider different classes of trees, where the class of a tree
is determined by the types and classes of her children. The
offspring distribution of a node depends on its type and on the
level of the tree where this node is living. In fact, we show that
the conditioned tree again is a multi-type Galton-Watson tree and
how this can be used to directly construct such a conditioned
tree. Our approach can be seen as a generalization of the
well-known decomposition of a supercritical Galton-Watson tree
into nodes whose offspring survives forever and nodes whose
offspring eventually goes to extinction as discussed in
\cite{LP:book}.

Section \ref{section:examples} discusses a couple of examples that
illustrate the applicability of our results. We give an example
concerning mutants in a population, we discuss an alternative to
Geiger's construction of a tree conditioned on having height
exactly $k$ and we show how to condition on the size of the $k$th
generation.

\subsection{Notation and preliminaries}\label{section:notation}

We will consider rooted multi-type Galton-Watson trees  with
arbitrary offspring distribution, that can depend on the current
generation. In such a tree each node has a type, which we indicate
by a natural number $t\in\Theta:=\left\{1,\ldots,\theta\right\}$.
If the root of a tree has type $t$, we use a bold-face $t$ to
denote this root. Define the set of trees of heigth $0$ as
$\mathcal{T}_0=\left\{{\mathbf 1},{\mathbf
2},\ldots,\boldsymbol{\theta}\right\}$. Then we define inductively
for $k\geq 1$ the set of trees of height at most $k$ by
\[ \T_k = \T_0\sqcup \bigsqcup_{n=1}^\infty \{(t,[T_1,\ldots, T_n])\mid T_i\in \T_{k-1},t\in\Theta\},\]
and denote the set of all trees by $\T = \bigcup_{k=0}^\infty
\T_k$. For a tree $T = (t,[T_1,\ldots, T_n]) \in
\Theta\times(\T_{k-1})^n$, the trees $T_1,\ldots, T_n$ will be
called the children of $T$ (notation: $T_i\prec T$). The type of
$T$ will be just the type of its root and will be denoted by
$r(T)$. We now define a function $N:\T\rightarrow\N^\theta$ that
counts how many children of each type a tree $T$ has:
$$
N(T) = (N_1(T),\ldots,N_\theta(T)),\quad\textrm{where}\quad N_t(T) = \#\left\{\tilde T\prec T:r(\tilde T)=t\right\}.
$$
The set of trees of heigth at most $k$ having a root of type $t$ is denoted by
$$
\T_{k}^t = \{T\in\T_k : r(T)=t\} = \{{\mathbf t}\}\sqcup \bigsqcup_{n=1}^\infty \{(t,[T_1,\ldots, T_n])\mid T_i\in \T_{k-1}\}.
$$
Let $\T^t = \bigcup_{k=0}^\infty \T_k^t$. Denote the offspring
distribution of a type $t$ node at height $l\geq 0$ by $\mu_l^t$,
for arbitrary probability measures $\mu_l^t$ on $\mathbb{N}^\theta
=\{0,1,\ldots\}^\theta$. Define independent random variables
$W_l^t\sim\mu_l^t$. For a vector ${\mathbf
x}\in\mathbb{N}^\theta$, we write the corresponding multinomial
coefficient as
$$
D({\mathbf x}) = \left(\begin{array}{c}|{\mathbf x}|_1\\
x_1,\ldots,x_\theta\end{array}\right) =
\frac{\left(\sum_{i=1}^\theta x_i\right)!}{\prod_{i=1}^\theta
x_i!}.
$$
We now introduce the Galton-Watson probability measures on
$\T_k^t$. Firstly, let $P_0^t$ be the trivial probability measure
on $\T_0^t$, so $P_0^t({\mathbf t}) = 1$. Now define inductively
the probability measure $P_{lk}^t$ for $0\leq l\leq k$ and $k\geq
1$ as the following probability measure on $\T_{k-l}^t$: if $l=k$,
then $P_{kk}^t=P_0^t$. Otherwise for all $T\in \T_{k-l}^t$
\[P_{lk}^t(T) = \frac{\P(W_{l}^t=N(T))}{D(N(T))}\prod_{\tilde T\prec T} P_{l+1,k}^{r(\tilde T)}(\tilde T),\]
where empty products are taken to be $1$. The intuition is that
the second sub-index determines the size of the final tree we are
considering, whereas the first sub-index determines at which level
we are building up the tree (so $P_{lk}^t$ generates trees of type
$t$ at level $l$ of size $k-l$). We are interested in $P_{0k}^t$,
which is the Galton-Watson probability measure on $\T_k^t$ (trees
cut off at height $k$ with a root of type $t$).\\

\noindent In the next section we will introduce a class of recursive-type events on which we would like to condition, and discuss several examples of such events. In Section \ref{section:conditionalmeasures} we will introduce the conditional measures corresponding to our events, and in Section \ref{section:proof} we show that these conditional measures indeed coincide with the original Galton-Watson measure, conditioned on our event.

\section{Conditioning on recursive events}\label{section:conditioning}

In this section we introduce a class of recursive-type events on which we would like to condition,
such as the event that the tree survives until a specific level.

\subsection{Partitioning the set of trees}

We will now set up our general framework and show how some
examples fit into it. We start by choosing $k_0\in\N$ and
partitioning $\T_{k_0}$ into $m$ classes
$A_{k_0}^{(1)},\ldots,A_{k_0}^{(m)}$. Typically, all trees in such
a class have some property that all trees in the other classes do
not have. One of the simplest examples would be a partition into
two classes, where trees that survive until some level $k$ are in
the first class and all other trees in the second. The partition
of $\T_{k_0}$ will be the starting point to recursively define
partitions of $\T_l$, $k_0<l\leq k$ into sets $A_l^{(i)},
i=1,\ldots,m$, where $k$ is the (maximum) height of the trees that
we are considering. Suppose that the partition of $\T_{l-1}$ is
already defined. Then we are able to introduce a counting matrix
for trees in $\T_l$. Define
$C_l:\T_l\rightarrow\N^{m\times\theta}$ such that for $T\in\T_l$
the $(i,j)$th position is given by
$$
C_l^{(i,j)}(T) = \#\left\{\tilde T\prec T:\tilde T\in A_{l-1}^{(i)}\cap \T^j\right\},
$$
so this is the number of children of $T$ having type $j$ and being
an element of the $i$th partition class. Now for $k_0<l\leq k$ we partition
$\N^{m\times\theta}$ into subsets $B_{l,1},\ldots,B_{l,m}$. This
partition is the key for the recursive definition of $A_l^{(i)}$.
The set $A_l^{(i)}$ will contain exactly those trees for which the
counting matrix $C_l(T)$ is in $B_{l,i}$:
\begin{equation}\label{eq:partition}
A_l^{(i)} = \left\{T\in\T_l:C_l(T)\in B_{l,i}\right\}.
\end{equation}

\subsection{Examples}\label{section:examples}

Before going into the details of the construction of conditioned
trees, we will discuss some examples of recursive events that can
be handled by our approach.

\subsubsection{Genetic mutations}
Suppose we have a population in which sometimes an individual
(mutant) is born having a particular mutation in its genetic
material. This mutation can be inherited by subsequent
generations. Suppose we know the probability that the root is
mutated. Such a population can be described as a two-type
Galton-Watson process in which the offspring distribution is type-
and possibly level-dependent. We take $\Theta =
\left\{1,2\right\}$ to be the set of types, where mutants have
type
$1$.\\
\\
Suppose we would like to condition on the event ``there is at
least one mutant in the $k$th generation". Choose $k_0=0$ and
partition $\T_0 = \left\{{\mathbf 1},{\mathbf 2}\right\}$ into the
classes $A_0^{(1)} = \left\{{\mathbf 1}\right\}$ and $A_0^{(2)} =
\left\{{\mathbf 2}\right\}$. For $0<l\leq k$, we want to define
$A_l^{(1)}$ and $A_l^{(2)}$ by
$$
A_l^{(1)} = \left\{T\in\T_l:\hbox{there is a mutant at level $l$}\right\},\qquad A_l^{(2)} = \T_l\backslash A_l^{(1)}.
$$
These events satisfy a recursive relation: $A_l^{(1)}$ contains
exactly those trees that have at least one child in
$A_{l-1}^{(1)}$. For $T\in \T_l$, the first row of the $2\times
2$-counting matrix $C_l(T)$ counts the  children of $T$ that are
in $A_{l-1}^{(1)}$. Therefore, for all $0<l\leq k$ we let
$$
B_{l,1}  = \left\{\left[\begin{array}{cc}a&b\\c&d\end{array}\right]\in \N^{2\times 2}:a+b\geq 1\right\}, B_{l,2} = \N^{2\times 2}\backslash B_{l,1},
$$
and now (\ref{eq:partition}) gives the desired partition of $T_l$.
For the sake of illustration, we note that with a minor change, we
can condition on ``there is at least one mutant in the $k$th
generation inheriting its mutation from the root". To achieve
this, it suffices to merely redefine $B_{l,1}$ and $B_{l,2}$ for
all $0<l\leq k$ as follows
$$
B_{l,1}  = \left\{\left[\begin{array}{cc}a&b\\c&d\end{array}\right]\in \N^{2\times 2}:a\geq 1\right\}, B_{l,2} = \N^{2\times 2}\backslash B_{l,1}.
$$

In these two examples, we defined one partition class $A_l^{(1)}\subseteq \T_l$ by the event on which conditioning is required. The only other partition class was just the complement of the first one. Finding a suitable partition of the set of trees is not always that obvious, as is demonstrated in the next example. We will show how to condition on the slightly more complicated event ``All mutants in the tree
inherit their mutation from the root and at least one mutant is
present in generation $k$". As before, define one partition class $A_l^{(1)}$ as the set of trees satisfying the condition. Here it is not sufficient to define only one other partition class. One obstacle is that some trees (namely those with a ``spontaneous mutation") in the complement $\T_l\setminus A_l^{(1)}$ are forbidden as a child of trees in $A_{l+1}^{(1)}$ and others are not.

Nevertheless, with a slightly more elaborate partition, we can still handle this case. We distinguish four
classes and partition $\T_0$ into
$$
A_0^{(1)} = \left\{{\mathbf 1}\right\},\quad A_0^{(2)} =
\left\{{\mathbf 2}\right\}\quad \textrm{and}\quad A_0^{(3)} =
A_0^{(4)} = \emptyset.
$$
For $0<l\leq k$, we define the following subsets of $\N^{4\times
2}$:
\begin{eqnarray*}
& B_{l,1}  =
\left\{\left[\begin{array}{cc}a&0\\b&c\\0&0\\d&0\end{array}\right]:a\geq
1\right\},&
 B_{l,2} = \left\{\left[\begin{array}{cc}0&0\\0&a\\0&0\\0&0\end{array}\right]\right\}\\
 &B_{l,3} =
\left\{\left[\begin{array}{cc}a&b\\c&d\\e&f\\g&h\end{array}\right]:b+e+f+h\geq
 1\right\},&
 B_{l,4} = \left\{\left[\begin{array}{cc}0&0\\a&b\\0&0\\c&0\end{array}\right]:a+c\geq 1\right\}\\
\end{eqnarray*}
As can be easily checked, these sets are disjoint and $\bigcup_{i}
B_{l,i} = \N^{4\times 2}$, so this indeed is a partition. It
follows by induction that the sets $A_l^{(i)}$ partition $\T_l$ in
such a way that
\begin{itemize}
\item $A_l^{(1)}, 0<l\leq k$ contains exactly the trees having
\begin{itemize}
\item at least one mutated child of which the mutated progeny reaches level $l$, and
\item no ``spontaneous" mutants in the progeny of their children.
\end{itemize}
\item $A_l^{(2)}, 0<l\leq k$ contains the trees having only type 2
descendants.

\item $A_l^{(3)}, 0<l\leq k$ contains the trees having a type 2
descendant with a type 1 child (``spontaneous mutation").

\item $A_l^{(4)}, 0<l\leq k$ contains all other trees in $\T_l$.
\end{itemize}
Note that these classes are defined by properties of the children of a
tree and not by the type of the tree itself. For example, a tree
in $A_l^{(2)}$ can have a type 1 root, but all its descendants
have type 2. The conditional measure we are interested in is now
obtained by conditioning $P_{0k}^1$ on $A_k^{(1)}$.

\subsubsection{Conditioning on the size of generation $k$}

As a next example, we show how to condition a single-type
Galton-Watson tree on having exactly $G$ individuals in the $k$th
generation. In this case, we partition $\T_0 = \left\{{\mathbf
1}\right\}$ into $G+2$ classes by defining
$$
A_0^{(1)} = \left\{{\mathbf 1}\right\},\quad A_0^{(0)} = A_0^{(2)} = A_0^{(3)} =\ldots =
A_0^{(G)} = A_0^{(G+1)} = \emptyset.
$$
Define $x\in\N^{G+2}$ by $x:=[0\quad 1\quad 2\quad\ldots\quad
G\quad G+1]^T$. For $0<\l\leq k$, we define
$$
B_{l,i} = \left\{y\in\N^{G+2}:x^Ty = i\right\},\quad B_{l,G+1} =
\left\{y\in\N^{G+2}:x^Ty \geq G+1\right\},
$$
where $0\leq i \leq G$. Partitioning $\T_l$ according to
(\ref{eq:partition}) gives the following: for $0\leq i\leq G$,
$A_l^{(i)}$ contains the trees of which the $l$th generation has
exactly size $i$, while $A_l^{(G+1)}$ contains the trees of which
the $l$th generation has at least size $G+1$. Conditioning
$P_{0k}$ on $A_k^{(G)}$ gives the result we are looking for.

\subsubsection{The tree has heigth exactly $k$}

As a final illustration, we explain how to condition a
Galton-Watson tree on having height exactly $k$, thus producing an
alternative for the construction of Geiger and Kersting
\cite{Geiger}. We consider trees in $\T_{k+1}$ that are
conditioned to reach level $k$, but not level $k+1$. We start by
choosing $k_0=2$, and partitioning $\T_2$ into three sets, namely
correct trees, short trees and long trees:
\begin{align*}
 A_2^{(1)} & =\{ T\in \T_2\mid T\mbox{ reaches level 1, but not level 2}\},\\
A_2^{(2)} & =\{T\in \T_2\mid T\mbox{ does not reach level 1}\} = \{\0\},\\
A_2^{(3)} & = \{T\in \T_2\mid T\mbox{ reaches level 2}\}.
\end{align*}
Define for each $2<l\leq k+1$
\begin{align*}
 B_{l,1} & = \{ n\in\N^3\mid n_1\geq 1, n_3=0\},\\
B_{l,2} & = \{ n\in\N^3\mid n_1=0, n_3=0\},\\
B_{l,3} & = \{ n\in\N^3\mid n_3\geq 1\},\\
\end{align*}
and let $\T_{l}$ be partitioned as in (\ref{eq:partition}). This
construction guarantees that if a tree $T\in \T_{k+1}$ is an
element of $A_{k+1}^{(1)}$, then it has at least one child that
reaches level $k$, and no children that reach level $k+1$. If
$T\in A_{k+1}^{(2)}$, all its children do not reach level $k$, and
if $T\in A_{k+1}^{(3)}$, then at least one child reaches level
$k+1$. Conditioning on being in $A_{k+1}^{(1)}$ therefore gives
the desired result.


\subsection{Remarks following the examples}

As it turns out from the examples in the previous section, the setup
allows to condition on quite a variety of events. A fundamental
requirement on these events is that they are determined only by
the \emph{number} of children of a tree having particular
properties. So we can (for instance) not distinguish between trees
having the same children in a different order.

An additional example is discussed in detail in \cite{Cator}. As an application of the theory developed in the present paper, the cost of searching a tree to a given
level is determined. The proposed model takes into account costs for having a lot of children, but also for
walking into dead ends. So both a high expected offspring and a low expected offspring would give high search costs. This gives rise to an optimization problem: which offspring distribution gives minimal costs? For this model the conditional probability measures are explicitly constructed, leading to recursions that
enable us to calculate the costs and solve the optimization problem for Poisson offspring.

Conditioning on recursive events as in the examples allows us to
compute (conditional) probabilities that are defined in terms of
such events. As an illustration: in the example on genetic
mutations we can easily compute the probability that the root is
mutated, given that there is at least one mutant in generation
$k$. What makes the results even more useful is that they show how
to directly construct a tree conditioned on some event. This means
that trees conditioned on (rare) events can be studied by just
simulating them.

\section{Conditional measures}\label{section:conditionalmeasures}

In this section we
construct an alternative measure $\tilde P_{lk}^t$ on
$\T_{k-l}^t$, that depends on the event we want to condition on.
As soon as we have this measure, conditioning on the desired event
is a triviality. In the next section, we will show that in fact
the two measures $P_{lk}^t$ and $\tilde P_{lk}^t$ are the same.\\

Define $^tp_{lk}^{(i)}$ for $0\leq l\leq k-k_0$ by
$$
^tp_{lk}^{(i)} = P_{lk}^t(A_{k-l}^{(i)}\cap \T^t).
$$
We can calculate this probability in a recursive way. Denote, for
$q\in [0,1]^m$ with $\sum q_i =1$, by ${\rm Multi}(n,q)\in\N^m$
the multinomial distribution where we distribute $n$ elements over
$m$ classes, according to the probabilities $q_i$. We also choose
independent random vectors $W_{l}^t\sim \mu_l^t$ according to the
offspring distribution of a type $t$ node at level $l$ and denote
the $j$th coordinate by $W_{l,j}^t$. Then, for $l<k-k_0$
\begin{equation}\label{eq:tplk(i)}
^tp_{lk}^{(i)} = \P\left( \bigotimes_{j=1}^\theta {\rm
Multi}\left(W_{l,j}^t,({^jp}_{l+1,k}^{(1)}, \ldots,
{^jp}_{l+1,k}^{(m)})\right) \in B_{k-l,i}\right),
\end{equation}
where, for
$a_0,\ldots,a_\theta\in\N^m$, we defined
$\Lambda:=\bigotimes_{j=1}^\theta a_j\in\N^{m\times\theta}$ to be
the matrix for which $\Lambda_{ij}= a_j(i)$.
%
%
%

We proceed by defining the conditional measure
$^t\tilde{Q}_{lk}^{(i)}$ on $A_{k-l}^{(i)}\cap \T^t$.
To do this, define for each $t\in\left\{1,\ldots,\theta\right\}$
and $0\leq l\leq k-k_0-1$ on the same probability space as $W_{l}^t$, the random matrices
$$
^tX_{lk} = \left(\begin{array}{ccc}
{^tX}_{lk}^{(1,1)} & \ldots & {^tX}_{lk}^{(1,\theta)}\\
\vdots & \ddots & \vdots \\
{^tX}_{lk}^{(m,1)} & \ldots & {^tX}_{lk}^{(m,\theta)}
\end{array}\right),
$$
such that conditional on $W_{l}^t$, all columns are independent
and the distribution of the $j$th column satisfies
\[ ({^tX}_{lk}^{(1,j)}, \ldots, {^tX}_{lk}^{(m,j)})\mid W_{l}^t\sim {\rm
Multi}\left(W_{l,j}^t,({^jp}_{l+1,k}^{(1)}, \ldots,
{^jp}_{l+1,k}^{(m)})\right).\] This determines the full joint
distribution of $(W^t_{l}, {^tX}_{lk})$. For a type $t$ node at
level $l$, the distribution of its children over the $\theta$
types is given by the random vector $W^t_{l}$. Furthermore, the
$j$th column of $^tX_{lk}$ represents how the type $j$ children of
this type $t$ node are distributed over the $m$ classes. For
$l=k-k_0$, we define for each $T\in A^{(i)}_{k_0}\cap \T^t$
\[ ^t\tilde{Q}_{k-k_0,k}^{(i)}(T) = \frac{P^t_{k-k_0,k}(T)}{P^t_{k-k_0,k}(A^{(i)}_{k_0}\cap \T^t)},\]
as a probability measure on $A^{(i)}_{k_0}\cap \T^t$. Next, we inductively define the probability measures $^t\tilde{Q}_{lk}^{(i)}$
on $A_{k-l}^{(i)}\cap \T^t$ for each $0\leq l\leq k-k_0-1$ such
that for each $T\in A^{(i)}_{k-l}\cap \T^t$
\[ ^t\tilde{Q}^{(i)}_{lk}(T) = \frac{\P({^tX}_{lk}= C_{k-l}(T)\mid {^tX}_{lk}\in B^{(i)}_{k-l})}{D(C_{k-l}(T))}\prod_{j=1}^m\prod_{\tilde{T}\prec T : \tilde{T}\in A_{k-l-1}^{(j)}} {^{r(\tilde T)}}\tilde{Q}^{(j)}_{l+1,k}(\tilde{T}),\]
where we extended the definition of $D$ (see Section
\ref{section:notation}) to integer-valued matrices, and once again
empty products are taken to be 1. Note that this definition is
valid for all $T\in\T_{k-l}^t$: we simply get
$^t\tilde{Q}_{lk}^{(i)}(T)=0$ whenever $T\not\in A_{k-l}^{(i)}$.
We can now define the alternative measure
$\tilde{P}^t_{lk}$ on $\T_{k-l}^t$:
\[ \tilde{P}^t_{lk}(T)= \sum_{i=1}^m {^tp}_{lk}^{(i)}{^t\tilde{Q}}_{lk}^{(i)}(T).\]

\subsection{Construction according to the conditional measure}

We can describe the random tree $\Tb\sim \tilde{P}_{lk}^t$ as
follows. The root of the tree has type $t$. To construct the tree,
we first toss an $m$-sided coin to determine in which of the $m$
classes $\Tb$ is, giving probability $^tp_{lk}^{(i)}$ to the $i$th
class $A_{k-l}^{(i)}$. If $\Tb\in A_{k-l}^{(i)}$, then we choose
it according to $^t\tilde{Q}_{lk}^{(i)}$. This means that we choose
$(\tilde{W}^t_{lk},{^t\tilde{X}}_{lk})$, where $\tilde{W}^t_{lk}$
counts the numbers of children of $\Tb$ of each type and $^t\tilde{X}_{lk}$
counts for each type the numbers of children that will lie in each
of the $m$ classes, according to
$$
(\tilde{W}^t_{lk},{^t\tilde{X}}_{lk})\sim
(W_{l}^t,{^tX}_{lk})\mid {^tX}_{lk}\in B^{(i)}_{k-l}.
$$
The $\sum_{i,j}{^t\tilde{X}}_{lk}^{(i,j)}$ children are
distributed over the $\sum_{j}\tilde{W}_{lk,j}^t$ positions
uniformly at random. Then for each child of type $j$ in
$A_{k-l-1}^{(i)}$ we draw a tree according to
$^j\tilde{Q}_{l+1,k}^{(i)}$.

In this way we have described the random tree as a Galton-Watson
tree with $m\times\theta$ `types' of children and type- and
level-dependent offspring distribution. Note that conditioning
$\tilde{P}_{lk}^t$ on $A_{k-l}^{(i)}$ is trivial: we simply have
to draw $\Tb$ according to $^t\tilde{Q}_{lk}^{(i)}$.

\section{The two random trees are equally distributed}\label{section:proof}

The following theorem shows that the construction procedure of
Section \ref{section:conditionalmeasures} in fact generates trees
with the same probabilities as under the original Galton-Watson
measure. Fix $k$ and $k_0$ and define all measures as before.
\begin{stelling} For all $0\leq l\leq k-k_0$, $t\in\Theta$ and $T\in \T_{k-l}^t$,
$$
P_{lk}^t(T)=\tilde{P}_{lk}^t(T).
$$
\end{stelling}
{\bf Proof:} The theorem is true by construction for $l=k-k_0$.
Now suppose that we have already shown that
$P_{l+1,k}^t=\tilde{P}_{l+1,k}^t$ for all $t$. Choose $T\in
\T_{k-l}$ and suppose $T\in A_{k-l}^{(i)}$. Before we show that
$\tilde P_{lk}^t(T)=P_{lk}^t(T)$, we collect some useful
observations. First of all, note that the number of ways to
distribute the individuals over the positions in $C_{k-l}$ can be
written as a product by first assigning a type to each individual
and then distributing all individuals of a given type over the
classes (writing $C_{k-l}$ for $C_{k-l}(T)$ and $N$ for $N(T)$):
$$
D(C_{k-l}) = D(N)\prod_{j=1}^\theta D\left(\left(C_{k-l}^{(i,j)}\right)_{i=1}^m\right).
$$
Secondly, note that $\P({^tX}_{lk}\in B_{k-l}^{(i)})$ is equal to
\begin{align*}
\sum_{(n_1,\ldots,n_\theta)\in\N^\theta}&\P\left(W_{l}^t =
(n_1,\ldots,n_\theta)\right)\P\left( \bigotimes_{j=1}^\theta {\rm
Multi}\left(n_j,({^jp}_{l+1,k}^{(1)}, \ldots,
{^jp}_{l+1,k}^{(m)})\right) \in B_{k-l,i}\right)\\
=\ &\P\left( \bigotimes_{j=1}^\theta {\rm
Multi}\left(W^t_{l,j},({^jp}_{l+1,k}^{(1)}, \ldots,
{^jp}_{l+1,k}^{(m)})\right) \in B_{k-l,i}\right)
\end{align*}
and by (\ref{eq:tplk(i)}) this is exactly ${^tp}_{lk}^{(i)}$. Next, since $C_{k-l}(T)$
determines $N(T)$, we have:
\begin{eqnarray*}
\P({^tX}_{lk} = C_{k-l}) &=& \P\left({^tX}_{lk}=C_{k-l}, W^t_{l}=N\right)\\
&=& \P\left(W^t_{l}=N\right)\P\left({^tX}_{lk}=C_{k-l}\mid W^t_{l}=N\right)\\
&=& \P\left(W^t_{l}=N\right)\prod_{j=1}^\theta \P\left(\left({^tX}_{lk}^{(i,j)}\right)_{i=1}^m = \left(C_{k-l}^{(i,j)}\right)_{i=1}^m\mid W^t_{l,j}=N_j\right)\\
&=& \P\left(W^t_{l}=N\right)\prod_{j=1}^\theta D\left(\left(C_{k-l}^{(i,j)}\right)_{i=1}^m\right) \prod_{i=1}^m \left({^jp}_{l+1,k}^{(i)}\right)^{C_{k-l}^{(i,j)}}\\
&=& \frac{\P\left(W^t_{l}=N\right)D(C_{k-l})}{D(N)}
\prod_{j=1}^\theta \prod_{i=1}^m
\left({^jp}_{l+1,k}^{(i)}\right)^{C_{k-l}^{(i,j)}}.
\end{eqnarray*}
Combining these observations gives
\begin{eqnarray*}
\tilde{P}_{lk}^t(T) &=& {^tp}_{lk}^{(i)}{^t\tilde Q}_{lk}^{(i)}(T)\\
&=& {^tp}_{lk}^{(i)}\frac{\P({^tX}_{lk}= C_{k-l}(T)\mid {^tX}_{lk}\in B^{(i)}_{k-l})}{D(C_{k-l}(T))}\prod_{j=1}^m\prod_{\tilde{T}\prec T : \tilde{T}\in A_{k-l-1}^{(j)}} {^{r(\tilde T)}}\tilde{Q}^{(j)}_{l+1,k}(\tilde{T})\\
&=& {^tp}_{lk}^{(i)}\frac{\P({^tX}_{lk}= C_{k-l}(T))}{D(C_{k-l}(T))\P({^tX}_{lk}\in B^{(i)}_{k-l})}\prod_{j=1}^m\prod_{\tilde t\in\Theta}\prod_{\tilde{T}\prec T : \tilde{T}\in A_{k-l-1}^{(j)}\cap \T^{\tilde t}} {^{\tilde t}}\tilde{Q}^{(j)}_{l+1,k}(\tilde{T})\\
&=& \left(\frac{\P\left(W^t_{l}=N\right)}{D(N)} \prod_{\tilde t=1}^\theta \prod_{j=1}^m \left({^{\tilde t}p}_{l+1,k}^{(j)}\right)^{C_{k-l}^{(j,\tilde t)}}\right)\left(\prod_{j=1}^m\prod_{\tilde t\in\Theta}\prod_{\tilde{T}\prec T : \tilde{T}\in A_{k-l-1}^{(j)}\cap \T^{\tilde t}} {^{\tilde t}}\tilde{Q}^{(j)}_{l+1,k}(\tilde{T})\right)\\
&=& \frac{\P\left(W^t_{l}=N\right)}{D(N)} \prod_{j=1}^m\prod_{\tilde t\in\Theta}\prod_{\tilde{T}\prec T : \tilde{T}\in A_{k-l-1}^{(j)}\cap \T^{\tilde t}} {^{\tilde t}p}_{l+1,k}^{(j)} {^{\tilde t}}\tilde{Q}^{(j)}_{l+1,k}(\tilde{T})\\
&=& \frac{\P\left(W^t_{l}=N\right)}{D(N)} \prod_{j=1}^m\prod_{\tilde t\in\Theta}\prod_{\tilde{T}\prec T : \tilde{T}\in A_{k-l-1}^{(j)}\cap \T^{\tilde t}} \tilde{P}_{l+1,k}^{\tilde t}(\tilde T)\\
&=& \frac{\P\left(W^t_{l}=N\right)}{D(N)} \prod_{\tilde{T}\prec T} P_{l+1,k}^{r(\tilde T)}(\tilde T)\\
&=& P_{lk}^t(T).
\end{eqnarray*}
\hfill $\Box$

\section{Example: genetic mutations revisited}

\begin{figure}
\begin{center}
\includegraphics[width = 10cm]{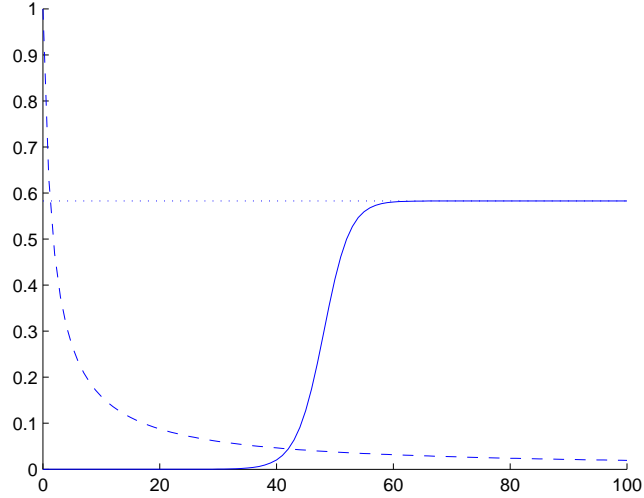}\caption{Dashed: ${^1}p_{0k}^{(1)}$. Solid: ${^2}p_{0k}^{(1)}$. Dotted: survival probability of a tree with Poisson($3/2$) offspring.}\label{fig:kansen}
\end{center}
\end{figure}

In this section we use our results to work out one of the examples
of Section \ref{section:examples}. For these calculations it will turn out to be very useful that our conditioned tree is again a Galton-Watson multitype tree. We consider a population with
mutants and let the set of types be $\Theta = \left\{1,2\right\}$,
where type $1$ denotes a mutant. The number of children of a type
$t$ node will have a ${\rm Pois}(\mu_t)$ distribution and each
child has probability $p_t$ to be a mutant itself, independent of
all other children.\\
\\
We will condition on the event that there is at least one mutant
in the $k$th generation. The corresponding partition of $\T_l$ is
given by
\begin{eqnarray*}
A_l^{(1)} &=& \left\{T\in\T_l:\hbox{there is a mutant at level
$l$}\right\} = \left\{T\in\T_l:C_l(T)\in B_{l,1}\right\},\\
A_l^{(2)} &=& \T_l\backslash A_l^{(1)},
\end{eqnarray*}
where $B_{l,1}$ is the set of matrices
$$
B_{l,1}  =
\left\{\left[\begin{array}{cc}a&b\\c&d\end{array}\right]\in
\N^{2\times 2}:a+b\geq 1\right\}.
$$
Remember that this means that a tree in $\T_l$ is an element of
$A_l^{(1)}$ if and only if it has a type $1$ or a type $2$ child
in $A_{l-1}^{(1)}$. We will now derive the recursions for the
probabilities $^tp_{lk}^{(i)}$. A type $t$ subtree that starts on
level $l<k$ has two types of children, and each child is in one of
the two classes. Type and class of a child are independent of all
other children's properties. So we can introduce four new `types',
occurring according to the following distributions:
\begin{equation}\label{eq:offspring}
\begin{array}{ll}
{^tX}_{lk}^{(1,1)} \sim \rm{Pois}\left({^1}p_{l+1,k}^{(1)}\cdot
p_t\cdot\mu_t\right), & {^tX}_{lk}^{(1,2)} \sim  \rm{Pois}\left({^2}p_{l+1,k}^{(1)}\cdot (1-p_t)\cdot\mu_t\right), \\
{^tX}_{lk}^{(2,1)} \sim\rm{Pois}\left({^1}p_{l+1,k}^{(2)}\cdot
p_t\cdot\mu_t\right),  & {^tX}_{lk}^{(2,2)} \sim
\rm{Pois}\left({^2}p_{l+1,k}^{(2)}\cdot (1-p_t)\cdot\mu_t\right),
\end{array}
\end{equation}
all independent of each other. In this notation
${^tX}_{lk}^{(a,b)}$ stands for the number of type $b$-children in
class $a$ of a type $t$ node at level $l$. Note that the
intensities indeed add up to $\mu_t$. 
For instance, the probability that a mutant (type $1$) on level
$l$ does not generate a mutant on level $k$ satisfies
$$
{^1}p_{lk}^{(2)} = e^{-{^1}p_{l+1,k}^{(1)}\cdot
p_1\cdot\mu_1}\cdot e^{-{^2}p_{l+1,k}^{(1)}\cdot
(1-p_1)\cdot\mu_1} = e^{-\mu_1\left({^1}p_{l+1,k}^{(1)}\cdot
p_1+{^2}p_{l+1,k}^{(1)}\cdot (1-p_1)\right)},
$$
and similarly it follows that
\begin{eqnarray*}
{^1}p_{lk}^{(1)} &=& 1-e^{-\mu_1\left({^1}p_{l+1,k}^{(1)}\cdot p_1+{^2}p_{l+1,k}^{(1)}\cdot (1-p_1)\right)},\\
{^2}p_{lk}^{(1)} &=& 1-e^{-\mu_2\left({^1}p_{l+1,k}^{(1)}\cdot p_2+{^2}p_{l+1,k}^{(1)}\cdot (1-p_2)\right)},\\
{^2}p_{lk}^{(2)} &=& e^{-\mu_2\left({^1}p_{l+1,k}^{(1)}\cdot
p_2+{^2}p_{l+1,k}^{(1)}\cdot (1-p_2)\right)}.
\end{eqnarray*}
The corresponding initial conditions are
$$
{^1}p_{kk}^{(1)} = {^2}p_{kk}^{(2)} = 1,\qquad {^1}p_{kk}^{(2)} =
{^2}p_{kk}^{(1)} = 0.
$$
Figure \ref{fig:kansen} shows the behavior of these probabilities for the following choice of parameters: mutants reproduce at rate $\mu_1 = 1$ and non-mutants at rate $\mu_2 = \frac{3}{2}$. Mutants can only generate mutants ($p_1 = 1$) and a child of a non-mutant has a very small probability to be a mutant, $p_2 = 10^{-9}$. The dashed line shows the probability that a tree with mutated root has a mutant on the $k$th level as a function of $k$. This is a critical tree with only mutants that eventually goes extinct. The solid line shows the probability that a tree with non-mutated root has a mutant on the $k$th level. This tree is supercritical, with reproduction rate (very close to) $\frac{3}{2}$. In a tree with Poisson($\mu$) offspring, the extinction probability is the non-trivial solution of
$$
s = \sum_{n=0}^\infty \frac{\mu^n}{n!}e^{-\mu}s^n = e^{\mu(s-1)}.
$$
For $\mu = \frac{3}{2}$, this gives $s\approx 0.417$. The corresponding survival probability is given as a dotted line in Figure \ref{fig:kansen}. This indicates that if the tree does not die out, then eventually there will be mutants almost surely, since the population grows exponentially. The population is of order $10^9$ around generation $\frac{9\log(10)}{\log(3/2)}\approx 51$, which explains the location of the increase of the solid line.

\begin{figure}
\begin{center}
\includegraphics[width = 10cm]{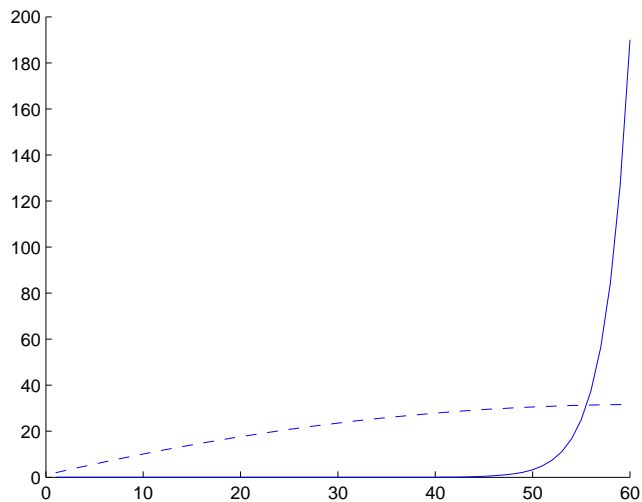}\caption{Expected number of mutants in each generation of a tree conditioned to have at least one mutant in generation $60$. Solid: the root is not a mutant. Dashed: the root is a mutant.}\label{fig:verwachtingen}
\end{center}
\end{figure}

We will now consider a tree with a non-mutated root, conditioned on having a mutant on
level $k$. We will show how to compute the expected number of
mutants on each level in such a conditioned tree. A type $t$
subtree at level $l$ is in class $1$ if and only if its counting
matrix is in $B^{(1)}_{k-l,1}$. This corresponds to conditioning
the distributions in (\ref{eq:offspring}) on
$$
{^tX}_{lk}^{(1,1)}+{^tX}_{lk}^{(1,2)} \geq 1.
$$
Write $X$ for ${^tX}_{lk}^{(1,1)}$ and $Y$ for
${^tX}_{lk}^{(1,2)}$. Then
$$
\E[X] = \E[X\mid X+Y \geq 1]\cdot\P(X+Y \geq 1) + \E[X\mid X+Y
=0]\cdot\P(X+Y =0).
$$
Since $\E[X\mid X+Y =0] = 0$ and $\P(X+Y\geq 1) =
{^t}p_{lk}^{(1)}$, we obtain
$$
\E[X\mid X+Y \geq 1] = \frac{\E[X]}{\P(X+Y \geq 1)} =
\frac{{^1}p_{l+1,k}^{(1)}\cdot p_t\cdot\mu_t}{^tp_{lk}^{(1)}}.
$$
And analogously:
$$
\E[Y\mid X+Y \geq 1] =  \frac{\E[Y]}{\P(X+Y \geq 1)} =
\frac{{^2}p_{l+1,k}^{(1)}\cdot (1-p_t)\cdot\mu_t}{^tp_{lk}^{(1)}}.
$$
The conditioned tree is in fact a four-type Galton-Watson tree:
\begin{eqnarray*}
\rm{type\ }\alpha: && \rm{mutant,\ class\ } 1\\
\rm{type\ }\beta: && \rm{mutant,\ class\ } 2\\
\rm{type\ }\gamma: && \rm{non-mutant,\ class\ } 1\\
\rm{type\ }\delta: && \rm{non-mutant,\ class\ } 2\\
\end{eqnarray*}
Abbreviating ${^a}p_{l+1,k}^{(b)}$ by ${^a}p^{b}$ and $(1-p_i)$ by
$q_i$, the expected offspring of a node at level $l<k$ described by
the following matrix:
$$ M_{lk} = \left(\begin{array}{cccc}
{^1}p^{1} \cdot p_1\cdot \mu_1/{^1p}_{lk}^{1} & 0 & {^1}p^{1}
\cdot p_2\cdot \mu_2/^2p_{lk}^{1} &
0\\
{^1}p^{2}\cdot p_1\cdot\mu_1 & {^1}p^{2}\cdot p_1\cdot\mu_1 &
{^1}p^{2}\cdot p_2\cdot\mu_2 &
{^1}p^{2}\cdot p_2\cdot\mu_2 \\
{^2}p^{1}\cdot q_1\cdot\mu_1/^1p_{lk}^{1} & 0 & {^2}p^{1}\cdot
q_2\cdot\mu_2/^2p_{lk}^{1} &
0\\
{^2}p^{2}\cdot q_1\cdot\mu_1 & {^2}p^{2}\cdot q_1\cdot\mu_1 &
{^2}p^{2}\cdot q_2\cdot\mu_2 &
{^2}p^{2}\cdot q_2\cdot\mu_2 \\
\end{array}\right).
$$
In this matrix, the columns give the expected offspring of a type
$\alpha$, $\beta$, $\gamma$ or $\delta$ node respectively at level
$l$ in the tree. Now the conditioned tree we are interested in is
just a tree with a root of type $\gamma$. The types $\alpha$ and
$\beta$ correspond to mutated individuals. Therefore, the expected
number of mutants on level $l$ in the conditioned tree is given by
$$
\left(1\quad 1\quad 0\quad 0\right) M_{l-1,k}M_{l-2,k}\ldots
M_{0,k}\left(\begin{array}{c}0\\0\\1\\0\end{array}\right).
$$
See Figure \ref{fig:verwachtingen}, for a plot of these expected numbers as a function of the generation.
\\

\begin{figure}
\begin{center}
\includegraphics[width = 11cm]{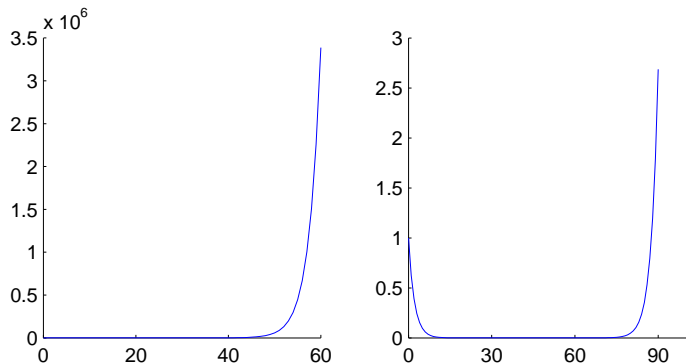}\caption{Expected size of each generation in a tree conditioned to have no mutant in generation $k$. Left: $k=60$, right: $k=90$. In both cases, the root is not a mutant. Note the huge difference in expected population size.}\label{fig:remarkable}
\end{center}
\end{figure}

As a last example, we computed the expected total number of individuals in a tree with a non-mutant root and conditioned to have \emph{no} mutant in generation $k$. If $k$ is small, occurrence of mutants is unlikely anyway, so then the tree just grows exponentially. If the population grows beyond order $10^9$, then the condition has a serious influence on the expected size of the tree. For example, taking $k=60$, the population size in the unconditioned tree would be of order $10^{10}$, but in the conditioned tree it is only of order $10^6$. See Figure \ref{fig:remarkable}, left plot. For even larger $k$ the condition to have no mutant in generation $k$ is very restrictive. See Figure \ref{fig:remarkable}, right plot. Apparently, the condition more or less forces the tree to die out early. In the unlikely case that it survives to generation $90$, the population stays small for a long time. The minimal expected size is $10^{-7}$ individuals and is attained around generation $40$. After that the population starts expanding. The later an individual is born, the less its progeny is influenced by the condition that no mutant is present in generation $90$. This also explains the increase in the curve that is seen at the end.

\section{Conclusion}

We have demonstrated how to condition multi-type Galton-Watson
trees on events having some recursive nature. More specifically,
we looked at partitions of the set of trees in which each
partition set is defined by some tree property. A crucial aspect
of these properties is that they are determined completely by
the types of children of the tree and the partition sets to which they
belong. As our examples show, there is a wide variety of events
fitting into this framework.

We have shown that such a conditioned tree itself is again a
multi-type Galton-Watson tree, and we derived equations for the
type- and level-dependent offspring distribution. These results
turn out to be very useful to analyze conditioned trees. Also, using our explicit construction procedure we can directly generate a tree that is
conditioned to satisfy some property that has very low probability,
which should also be useful for simulation purposes.

\end{document}